\documentclass[a4paper, 9pt]{article}
\usepackage{latexsym}
\usepackage{amsmath}
\usepackage{graphics}
\usepackage{amsthm}
\usepackage[vcentermath]{youngtab}
\usepackage{amssymb}
\usepackage[all]{xy}
\newfont{\Fr}{eufm10}
\newfont{\Sc}{eusm10}
\newfont{\Bb}{msbm10}
\newfont{\Am}{msam10}
\newfont{\am}{msam7}
\numberwithin{equation}{section}
\newtheorem{theorem}{Theorem}[section]
\newtheorem{proposition}[theorem]{Proposition}
\newtheorem{lemma}[theorem]{Lemma}
\newtheorem{corollary}[theorem]{Corollary}
\newtheorem{claim}{Claim}

\newtheorem{ftheorem}{Theorem}{\bf}{\it}
{\bf}{\it}
\theoremstyle{definition}
\newtheorem{definition}[theorem]{Definition}

{\bf}{\rm}
\newtheorem{definition and corollary}[theorem]{Definition and Corollary}
\newtheorem{definition and lemma}[theorem]{Definition and Lemma}

\theoremstyle{remark}
\newtheorem{example}[theorem]{Example}
\newtheorem{remark}[theorem]{Remark}

\newtheorem{fexample}[ftheorem]{Example}{\it}{\rm}

\newcommand{\Lie}{\mbox{\rm Lie }}

\newcommand{\h}{\mathfrac{\h}}

\title{An exotic Springer correspondence for symplectic groups\footnote{This is a revised version of a preprint entitled ``On the geometry of exotic nilpotent cones", available through math.RT/0607478v1. We changed the title in order to stress our central acheivement.}}
\author{Syu \textsc{Kato}
          \footnote{Graduate School of Mathematical Sciences, University of Tokyo, 3-8-1 Meguro Komaba 153-8914, Japan.} \footnote{The author was supported by the JSPS Research Fellowship for young scientists during this research.}}

\begin{document}
\maketitle

\begin{abstract}
Let $G$ be a complex symplectic group. In [K1], we singled out the nilpotent cone $\mathfrak N$ of some reducible $G$-module, which we call the ($1$-) exotic nilpotent cone. In this paper, we study the set of $G$-orbits of the variety $\mathfrak N$. It turns out that the variety $\mathfrak N$ gives a variant of the Springer correspondence for Weyl groups of type $C$, but shares a similar flavor with that of type $A$ case. (I.e. there appears no non-trivial local system and the correspondence is bijective.) As an application, we present one sufficient condition for the bijectivity of our exotic Deligne-Langlands correspondence [K1].
\end{abstract}

\section{Main results}
Let $G = \mathop{Sp} ( 2n, \mathbb C )$ be a complex symplectic group. Let $B$ and $T$ be its Borel subgroup and a maximal torus of $B$, respectively. We denote by $X ^* ( T )$ the character group of $T$. Let $R$ be the root system of $(G, T)$ and let $R ^+$ be its positive part defined by $B$. Let $W := N _G ( T ) / T$ be the Weyl group of $( G, T )$. For a group or an algebra $H$, we put $\mathsf{Irr} H$ the set of isomorphism classes of simple $H$-modules. We embed $R$ and $R ^+$ into a $n$-dimensional Euclid space $\mathbb E = \oplus _i \mathbb C \epsilon _i$ as:
$$R ^+ = \{ \epsilon _i \pm \epsilon _j \} _{i < j} \cup \{ 2 \epsilon _i \} \subset \{ \pm \epsilon _i \pm \epsilon _j \} \cup \{ \pm 2 \epsilon _i \} = R \subset \mathbb E.$$
We define $V _1 := \mathbb C ^{2n}$ and $V _2 := ( \wedge ^2 V _1 ) / \mathbb C$. These representations have $B$-highest weights $\epsilon _1$ and $\epsilon _1 + \epsilon _2$, respectively. We put $\mathbb V := V _1 \oplus V _2$. For a $G$-module $V$, we define its weight $\lambda$-part (with respect to $T$) as $V [ \lambda ]$. The positive part $V ^+$ of $V$ is defined as
$$V ^+ := \bigoplus _{\lambda \in \mathbb Q _{\ge 0} R ^+ - \{ 0 \} } V [ \lambda ].$$
We define
$$F := G \times ^B \mathbb V ^+ \subset G \times ^B \mathbb V \cong G / B \times \mathbb V.$$
Composing with the second projection, we have a map
$$\mu : F \longrightarrow \mathbb V.$$
We denote the image of $\mu$ by $\mathfrak N$. This is the $G$-variety which we refer as the {\it exotic nilpotent cone}\footnote{From the view-point of invariant theory, our variety $\mathfrak N$ is nothing special. (See eg. Schwarz \cite{Sc} and Dadok-Kac \cite{DK}) The author considers it as an ``exotic nilpotent cone of the Lie algebra $\mathfrak g$" since it shares many representation-theoretic features with the nilpotent cone of $\mathfrak g$ in an ``exotic" fashion. Theorem C of this paper represents one of such features.}. By abuse of notation, we may denote the map $F \rightarrow \mathfrak N$ also by $\mu$. Basic properties of $\mathfrak N$ are:

\begin{ftheorem}[Geometric properties of $\mathfrak N$, cf. \cite{K} 1.2]\label{fgeom}
We have:
\begin{enumerate}
\item The $G \times \mathbb ( \mathbb C ^{\times} ) ^2$-action on $\mathbb V$ descends to $\mathfrak N$;
\item The defining ideal of $\mathfrak N$ is $( \mathbb C [ \mathbb V ] ^G _+ ) \mathbb C [ \mathbb V ] = ( \mathbb C [ V _2 ] ^G _+ ) \mathbb C [ \mathbb V ]$;
\item The variety $\mathfrak N$ is normal;
\item The map $\mu$ is a birational projective morphism;
\item The variety $Z := F \times _{\mathfrak N} F$ is a union of equi-dimensional irreducible varieties.
\end{enumerate}
\end{ftheorem}

This result itself follows rather easily from the results and technique developed by Kempf, Schwarz, Hesselink, and Dadok-Kac. Here we repeat the statement for the convenience of readers.

Let $\mathcal P _2 ( n )$ denote the set of pairs $( \lambda, \mu )$ of partitions such that $\left| \lambda \right| + \left| \mu \right| = n$. It is well-known that the set $\mathcal P _2 ( n )$ parametrizes $\mathsf{Irr} W$ (cf. Macdonald \cite{Mc} I Appendix B).

\begin{ftheorem}[= Corollary \ref{cormain}]\label{fJNF}
The set of $G$-orbits of $\mathfrak N$ is in one-to-one correspondence with $\mathcal{P} _2 ( n )$.
\end{ftheorem}

The proof is divided into three steps: First, we introduce a set $\mathcal{MP} ( n )$ (cf. 2.2) and construct a bijection with $\mathcal P _2 ( n )$. This is done via an intermediate set $\mathcal{SP} ( n )$ (cf. 3.3). Then, we construct a map $\mathcal{MP} ( n ) \mapsto G \backslash \mathfrak N$ and show it is surjective. This surjectivity statement is an enhancement of a correspondence given by Sekiguchi-Ohta \cite{Se, Oh}. Finally, we deduce the injectivity as a byproduct of Theorem \ref{Spr} explained in the below.

Let $\mathbb O$ be a $G$-orbit in $\mathfrak N$. Let $X \in \mathbb O$. By means of the Ginzburg theory \cite{CG, Gi} and \cite{K} 2.13, we equip the vector space
$$M _{\mathbb O} := H _{\mathrm{codim} \mathbb O} ( \mu ^{-1} ( X ), \mathbb C )$$
with an action of $W$. (The symbol $H _{\bullet}$ always represents the Borel-Moore homology group instead of the usual homology group.) Here the RHS is independent of the choice of $X \in \mathbb O$ (as $W$-modules).

\begin{ftheorem}[= Theorem \ref{main}]\label{Spr}
The assignment $\mathbb O \mapsto M _{\mathbb O}$ establishes a one-to-one correspondence between the set of $G$-orbits of $\mathfrak N$ and  $\mathsf{Irr} W$.
\end{ftheorem}

This is an enhancement of Grinberg's generalized Springer correspondence \cite{Gr}. Apart from the Ginzburg theory and a weak form of Theorem \ref{fJNF}, the proof depends on two facts: One is the numerical identity $\left| \mathcal{MP} ( n ) \right| = \left| \mathsf{Irr} W \right|$, and the other is the connectedness of the stabilizer of $G$-orbits of $\mathfrak N$. The former is a numerical version of  a combinatorial bijection $\mathcal{P} _2 ( n ) \leftrightarrow \mathcal{MP} ( n )$ defined as above. The latter follows by a result of Igusa (cf. Theorem \ref{Ig}).

For each $\lambda \in X ^* ( T ) \backslash \{ 0 \}$, we fix a basis element $\mathbf v [ \lambda ] \in \mathbb V [ \lambda ]$. Then, the most simple non-trivial example of Theorem \ref{Spr} is:

\begin{fexample}[$n = 2$] Let $G = \mathop{Sp} ( 4, \mathbb C )$ and let $\mathcal N$ be the nilpotent cone of $\mathfrak{sp} ( 4, \mathbb C )$. Let $\mathbf x [ \alpha ] \in \mathfrak{sp} ( 4, \mathbb C )$ be a non-zero $T$-root vector corresponding to $\alpha$. We put $\alpha _1 := \epsilon _1 - \epsilon _2$. We have:\\
\begin{center}
\begin{tabular}{ccc|c|c}
$\mathsf{Irr} W$ & $G \backslash \mathfrak N$ & $\mathcal P _2 ( 2 )$ & dim & $G \backslash \mathcal N$\\ \hline
sign & \{0\} & $( {\tiny \young(\empty,\empty)}, \emptyset )$ & 1 & \{0\}\\
Ssign & $\mathbf v [ \epsilon _ 1 ]$ & $( \emptyset, {\tiny \young(\empty,\empty)} )$ & 1 & $\mathbf x [ 2 \epsilon _1 ]$\\
Lsign & $\mathbf v [ \alpha _1 ]$ & $( {\tiny \yng(2)}, \emptyset )$ & 1 & $\mathbf x [ \alpha _1 ]$\\
regular & $\mathbf v [ \alpha _1 ] + \mathbf v [ \epsilon _ 1 ]$ & $( {\tiny \young(\empty)}, {\tiny \young(\empty)} )$ & 2 & $\mathbf x [ \alpha _1 ]$ \\
triv & $\mathbf v [ \alpha _1 ] + \mathbf v [ \epsilon _ 2 ]$ & $( \emptyset, {\tiny \yng(2)} )$ & 1 & $\mathbf x [ \alpha _1 ] + \mathbf x [ 2 \epsilon _2 ]$\\\hline
\end{tabular}
\end{center}
Here the sets $G \backslash \mathfrak N$ and $G \backslash \mathcal N$ specify the corresponding orbits in our result and the Springer correspondence, respectively. If we denote the usual Springer resolution by $\mu _S$, the above table implies:
$$H _2 ( \mu ^{-1} _S ( \mathbf x [ \alpha _1 ] ) ) \cong \text{Lsign} \oplus \text{regular}$$
and
$$H _4 ( \mu ^{-1} ( \mathbf v [ \alpha _1 ] ) ) \cong \text{Lsign}, H _2 ( \mu ^{-1} ( \mathbf v [ \alpha _1 ] + \mathbf v [ \epsilon _ 1 ] ) ) \cong \text{regular}.$$
\end{fexample}

We set $Z := F \times _{\mathfrak N} F$. For $a = ( s, q _1, q _2 ) \in G \times ( \mathbb C ^{\times} ) ^2$, we define $\mathfrak N ^a$ and $Z ^a$ to be the subvarieties of $\mathfrak N$ and $Z$ consisting of $a$-fixed points, respectively. We set $G ( s ) := Z _G ( s )$. With an aid of Theorem \ref{Spr}, we prove:

\begin{ftheorem}[= Theorem \ref{reg}]\label{cDL}
Let $a = ( s, q _1, q _2 ) \in G \times ( \mathbb C ^{\times} ) ^2$ be a semi-simple element such that $q _2$ is not a root of unity. Then, there exists a one-to-one correspondence
$$G ( s ) \backslash \mathfrak N ^a \leftrightarrow \mathsf{Irr} H _{\bullet} ( Z ^a ),$$
where $H _{\bullet} ( Z ^a )$ acquires an associative algebra structure by means of convolution operations $($cf. \S 4 or $\cite{CG}$ \S 2$)$.
\end{ftheorem}

In \cite{K}, we study the representation theory of the affine Hecke algebra $\mathbb H$ of type $C _n$. It is an associative algebra with three independent parameters $( q _0, q _1, q _2)$. Its quotient by the two-sided ideal generated by $( q _0 + q _1 )$ is isomorphic to the extended affine Hecke algebras of type $B _n$. In the language of \cite{K}, Theorem \ref{cDL} implies the following:

\begin{ftheorem}[= Corollary \ref{BDL}]
For the extended affine Hecke algebra of type $B _n$ with two-parameters $( - q _1, q _1, q _2 )$, the regularity condition of parameters holds automatically unless $- q _1 ^2 \neq q _2 ^{\pm m}$ for $0 \le m < n$ or $q _2 ^l \neq 1$ for $1 \le l < 2n$.
\end{ftheorem}

A large part of the proofs of the above two theorems are borrowed from Lusztig \cite{L3} with several minor modifications. Most notably, we do not need the odd-term vanishing condition of the homology of fibers by using Borho-MacPherson's argument and localization technique. We do not know whether the odd-term vanishing result holds in this setting. (But it is expected to have a close connection to the project of Achar-Henderson \cite{AH}.)

In order to deepen the subjects of \cite{K} and this paper, we need to determine the correspondence of Theorem \ref{Spr} explicitly. This is done in a subsequent paper \cite{K2}, where we employ a different method as well as the results of this paper.

\section{Preparatory materials}
In this section, we collect some extra notation which we use throughout this paper. For each $X \in \mathbb V$, we write
$$X := \mathbf v [ 0 ] +\sum _{\lambda \in X ^* ( T )} X  ( \lambda ) \mathbf v [\lambda],$$
where $\mathbf v [ 0 ] \in \mathbb V [ 0 ]$. We define the support of $X$ as
$$\| X \| := \{ i \in [ 1, n ] : X ( \pm \epsilon _i ) \neq 0 \text{ or } X ( \pm \epsilon _i \pm \epsilon _j ) \neq 0 \text{ for some sign and } j \in [1, n] \}.$$

The following definition of normal form is a slight enhancement of the good basis of Ohta \cite{Oh}. (See also Igusa \cite{Ig}.)

\begin{definition}[Normal forms]
A block $J ( \lambda )$ of length $\lambda$ and position $i$ is one of the following vectors in $\mathbb V ^+$:
$$\mathbf v ^{(j)} ( \lambda ) _{i} := ( 1 - \delta _{j, 0} ) \mathbf v [ \epsilon _{i + j} ] + \sum _{k = 1} ^{\lambda - 1} \mathbf v [ \alpha _{i + k}],$$
where $\alpha _i := \epsilon _i - \epsilon _{i + 1}$, $0 \le j \le \lambda$ is an integer, and $\delta _{j, 0}$ is Kronecker's delta. It is clear that $\| \mathbf v ^{(j)} ( \lambda ) _{i} \| = [i + 1, \lambda + i ]$ or $\emptyset$. A normal form of $\mathbb V$ is a sum
$$v = \sum _{i = 0} ^{n - 1} \mathbf v ^{(j _{i})} ( \lambda _{i} ) _{i} \in \mathbb V ^+$$
such that $\| \mathbf v ^{(j _{i})} ( \lambda _{i}  ) _{i} \| \cap \| \mathbf v ^{(j _{i ^{\prime}})} ( \lambda _{i ^{\prime}} ) _{i ^{\prime}} \| = \emptyset$ if $i \neq i ^{\prime}$.
\end{definition}

\begin{definition}[Marked partitions]
A marked partition $\vec{\lambda} = ( \lambda, a )$ is a partition $\lambda = ( \lambda _1 \ge \lambda _2 \ge \ldots )$ of $n$, together with a sequence $a = ( a _1, a _2, \ldots )$ of integers such that:
\begin{enumerate}
\item $0 \le a _k \le \lambda _k$ for each $k$;
\item $a _k = 0$ if $\lambda _{k + 1} = \lambda _k$;
\item $\lambda _p - \lambda _q > a _p - a _q > 0$ if $p < q$ and $a _p \neq 0 \neq a _q$.
\end{enumerate}
We denote the set of marked partitions by $\mathcal{MP} ( n )$. For $( \lambda, a ) \in \mathcal{MP} ( n )$, we put
$$J ( \lambda, a ) := \sum _{p \ge 1} \mathbf v ^{( a _p )} ( \lambda _p ) _{\lambda _{p} ^{<}},$$
where $\lambda _{p} ^{<} = \sum _{q < p} \lambda _q$. It is clearly a normal form.
\end{definition}

For a $G$-variety $\mathcal X$, we denote the set of $G$-orbits of $\mathcal X$ by $\mathfrak O _{\mathcal X}$.

\begin{remark}[Ohta \cite{Oh} \S 2]
Let $\mathbf 0$ denote the sequence $( 0, 0, \ldots )$. Then, the assignment
$$\mathcal{P} ( n ) \ni \lambda \mapsto G . J ( \lambda, \mathbf 0 ) \in \mathfrak O _{\mathfrak N \cap V _2 }$$
gives a one-to-one correspondence.
\end{remark}


\section{Combinatorial correspondence}
We retain the setting of the previous section. In this section, we present combinatorics which is needed in the sequel.

\begin{definition}[Bi-partition]
Let $\mathcal P _2 ( n )$ denote the set of pairs $( \lambda, \mu )$ of partitions such that $\left| \lambda \right| + \left| \mu \right| = n$.
\end{definition}

It is known that $\mathsf{Irr} W$ is parametrized by $\mathcal P _2 ( n )$ (cf. \cite{Mc} I Appendix B).

\begin{theorem}\label{Comb}
The set $\mathcal P _2 ( n )$ is in one-to-one correspondence with $\mathcal{MP} ( n )$ 
\end{theorem}

\begin{corollary}
The set $\mathcal{MP} ( n )$ is in one-to-one correspondence with $\mathsf{Irr} W$. \hfill $\Box$
\end{corollary}

The rest of this section is devoted to the construction of a bijection between $\mathcal{MP} ( n )$ and $\mathcal P _2 ( n )$. For two partitions $\lambda, \mu$, we define their sum $\lambda \odot \mu$ as the partition $\{ \lambda _p + \mu _p \} _{p \ge 1}$.

\begin{definition}
Let $I \subset \mathbb Z$. A sub-segment $[r, s] \subset I$ is called a component if $r - 1 \not\in I \not\ni s + 1$. A pair $( \lambda, I )$ is called a segmented partition of $n$ if $\lambda$ ($= ( \lambda _1 \ge \lambda _2 \ldots )$) is a partition of $n$ and $I \subset [ 1, \lambda _1 ]$ (possibly empty) decomposes into components
$$I = [ i _1, \lambda _{j _1} ] \cup [ i _2, \lambda _{j _2}] \cup \cdots$$
such that $i _p < \lambda _{j _p} + 1 < i _{p + 1}$ holds for every possible $p$. We denote the set of segmented partitions by $\mathcal{SP} ( n )$.
\end{definition}

\begin{example}[$n = 3$]We have the following correspondences:
\begin{center}
\begin{tabular}{c|cccc|cc}
$\lambda$ & (3) & (3) & (3) & (3) & (1, 1, 1) & (1, 1, 1)\\ \hline
$\mathcal{MP} ( 3 )$ & $\mathbf 0$ & $( 1, 0 )$ & $(2, 0)$ & $(3, 0)$ & $\mathbf 0$ & $(0,0,1)$ \\
$\mathcal{SP} ( 3 )$ & $\emptyset$ & $[1, 3]$ & $[2, 3]$ & $[ 3, 3 ]$ & $\emptyset$ & $[1, 1]$\\
$\mathcal{P} _2 ( 3 )$ & $\{ ( 3 ), \emptyset \}$ & $\{ \emptyset, (3) \}$ & $\{(1), (2) \}$ & $\{ (2), (1)\}$ & $\{ (1, 1, 1), \emptyset \}$& $\{\emptyset, (1,1,1)\}$\\ \hline
\end{tabular}
\end{center}
\begin{center}
\begin{tabular}{c|cccc}
$\lambda$ & (2, 1) & (2, 1) & (2, 1) & (2, 1) \\ \hline
$\mathcal{MP} ( 3 )$ & $\mathbf 0$ & $( 1, 0 )$ & $(2, 0)$ & $(0, 1)$ \\
$\mathcal{SP} ( 3 )$ & $\emptyset$ & $[1, 2]$ & $[2, 2]$ & $[ 1, 1 ]$\\
$\mathcal{P} _2 ( 3 )$ & $\{ (2, 1), \emptyset \}$ & $\{ \emptyset, (2, 1)\}$ & $\{(1, 1), ( 1 )\}$ & $\{(1), (1, 1)\}$\\ \hline
\end{tabular}
\end{center}

\end{example}

\begin{lemma}\label{SPMP}
There exists a one-to-one correspondence between the sets $\mathcal{MP} ( n )$ and $\mathcal{SP} ( n )$.
\end{lemma}

\begin{proof}
Let $\lambda$ be a partition of $n$. Let $\mathbf{x} = ( x _1, x _2, \ldots )$ be a strictly decreasing finite positive integer sequence. We define
$$\mathcal{MP} ( \lambda, \mathbf{x} ) := \{ ( \lambda, a ) \in \mathcal{MP} ( n ) : a _p \neq 0 \Leftrightarrow p \in \mathbf{x} \}.$$
Notice that $\mathcal{MP} ( \lambda, \mathbf{x} ) = \emptyset$ if $\mathbf{x}$ contains $p$ such that $\lambda _p = \lambda _{p + 1}$. Let $\mathcal{SP} ( \lambda, \mathbf{x} )$ be the set of segmented partitions such that $\{ \lambda _p \} _{p \in \mathbf x}$ is the set of (right) boundaries of components. In other word, we define
$$\mathcal{SP} ( \lambda, \mathbf{x} ) := \{ ( \lambda, I ) \in \mathcal{SP} ( n ) : \mathbf x = \{ p \in \mathbb Z : \lambda _p \in I \not\ni \lambda _p + 1, \lambda _{p + 1} \neq \lambda _p \} \}.$$
We have
$$\mathcal{MP} ( n ) = \bigsqcup _{\lambda \in \mathcal P ( n ), \mathbf{x}} \mathcal{MP} ( \lambda, \mathbf{x} ) \text{ and } \mathcal{SP} ( n ) = \bigsqcup  _{\lambda \in \mathcal P ( n ), \mathbf{x}} \mathcal{SP} ( \lambda, \mathbf{x} ).$$
\begin{claim}
We have $\# \mathcal{MP} ( \lambda, \mathbf{x} ) = \lambda _{x _1} \prod _{i = 2} ^{\# \mathbf{x}} (\lambda _{x _i} - \lambda _{x _{i + 1}} - 1 )$.
\end{claim}

\begin{proof}
We count the number of possible markings $a = ( a _1, a _2, \ldots )$. Since we have $a _p = 0$ if $p \not\in \mathbf{x}$, we restrict our attension to the set $\{ a _p \} _{p \in \mathbf{x}}$. We count the possible choice of $a _{x _i}$ with the knowledge of $a _{x _1}, a _{x _2}, \ldots, a _{x _{i - 1}}$. If $i = 1$, then the possible choice is $1 \le a _{x _1} \le \lambda _{x _1}$ since $a _{x _1} \neq 0$. If $i > 1$, then the possible choice is
$$\max \{ a _{x _j} : j < i \} < a _{x _i} < \min \{ \lambda _{x _i} - \lambda _{x _j} + a _{x _j} : j < i \}.$$
By definition, both the maximal and the minimal are attained at $j = i - 1$. Therefore, the number of possible choice of $a _i$ is $(\lambda _{x _{i - 1}} - \lambda _{x _{i}} - 1 )$ (independent of $\{ a _j \} _{j > i}$). Therefore, multiplying these yields the result.
\end{proof}

\begin{claim}
We have $\# \mathcal{SP} ( \lambda, \mathbf{x}) = \lambda _{x _1} \prod _{i = 2} ^{\# \mathbf{x}} (\lambda _{x _i} - \lambda _{x _{i + 1}} - 1 )$.
\end{claim}

\begin{proof}
By definition, we have a non-empty segment for each $\{ \lambda _{q} \} _{q \in \mathbf{x}}$. By definition, the number of choice of segments ended at $\lambda _{x _{i}}$ is $\lambda _{x _i}$ ($i = 1$) or $(\lambda _{x _i} - \lambda _{x _{i + 1}} - 1 )$ ($i \ge 2$). Multiplying these yields the result.
\end{proof}

By the comparison of two claims, we deduce
$$\# \mathcal{MP} ( \lambda, \mathbf{x} ) = \# \mathcal{SP} ( \lambda, \mathbf{x}).$$
Hence, it suffices to construct an injective assignment $\mathcal{SP} ( \lambda, \mathbf{x}) \rightarrow \mathcal{MP} ( \lambda, \mathbf{x} )$. Fix $( \lambda, I ) \in \mathcal{SP} ( \lambda, \mathbf{x})$. We construct an integer sequence $a ( \lambda, I ) = ( a _1, a _2, \ldots )$ inductively as follows:
\begin{enumerate}
\item If $p > n$, then we put $a _p := 0$;
\item Assume that $\{ a _{q} \} _{q > p}$ is already defined. Then, we put
$$a _p := \begin{cases} \min \{ j, j - \lambda _q + a _q - 1 : q > p, a _q \neq 0\}& \text{(if $[ j, \lambda _p ] \subset I$ is a component)} \\ 0 & \text{(otherwise)} \end{cases}$$
\end{enumerate}
We have
$$a _p \ge ( \lambda _q  + 2 ) - \lambda _q + a _q - 1 > a _q \text{ if $p < q$ and $a _p \neq 0 \neq a _q$}.$$
Moreover, we deduce
$$a _p \le j - \lambda _q + a _q - 1 < \lambda _p - \lambda _q + a _q$$
if $[ j, \lambda _p ] \subset I$ is a component. In particular,
$$\lambda _p - a _p > \lambda _q - a _q \text{ if $p < q$ and $a _p \neq 0 \neq a _q$}.$$
Hence, we deduce $( \lambda, a ( \lambda, I ) ) \in \mathcal{MP} ( \lambda, \mathbf{x} )$. Thus, it suffices to show $a ( \lambda, I ) = a ( \lambda, I ^{\prime} )$ only if $I = I ^{\prime}$. Assume that $I \neq I ^{\prime}$. Then, there exists maximal number $q$ such that 1) $a _p = a _{p} ^{\prime}$ for every $p > q$ and 2) $[ j, \lambda _q ]$ is a component of $I$ and not a component of $I ^{\prime}$. Then, we have necessarily $a _q \neq a _q ^{\prime}$. Therefore, the assignment
$$\mathcal{SP} ( n ) \ni ( \lambda, \mathbf{x} ) \mapsto ( \lambda, a ( \lambda, I ) ) \in \mathcal{MP} ( \lambda, \mathbf{x} )$$
determines an injective map as desired. 
\end{proof}

\begin{lemma}\label{P2SP}
There exists a one-to-one correspondence between the sets $\mathcal{P} _2 ( n )$ and $\mathcal{SP} ( n )$.
\end{lemma}

\begin{proof}
First, we construct a map $\mathcal P _2 ( n ) \rightarrow \mathcal{SP} ( n )$. Fix $( \lambda, \mu ) \in \mathcal{P} _2 ( n )$. We define $j _p ^+ := \lambda _p + \mu _p$ and $j _p := \lambda _p + \mu _{p + 1} + 1$. Then, we define
$$I ( \lambda, \mu ) := \bigcup _{p \ge 1} [ j _p, j _p ^+ ].$$
Here we understand $[ j _p, j _p ^+ ] = \emptyset$ when $j _p ^+ < j _p$. We regard $[ a, b ] \cup [b + 1, c] = [ a, c ]$ as segments. We have $j _p ^+ \in \lambda \odot \mu$. Thus, we conclude $( \lambda \odot \mu, I ( \lambda, \mu )  ) \in \mathcal{SP} ( n )$.\\
Second, we construct a map $\mathcal{SP} ( n ) \rightarrow \mathcal P _2 ( n )$. Fix $( \lambda, I ) \in \mathcal{SP} ( n )$. Put
$$\mu _p ( \lambda, I ) := \# \left( I \cap [ 1, \lambda _p ] \right ), \gamma ( \lambda, I ) := \lambda _p - \mu _p ( \lambda, I ).$$
It is clear that $( \mu ( \lambda, I ), \gamma ( \lambda, I ) ) \in \mathcal{P} _2 ( n )$. These two maps are mutually inverse, which implies the result.
\end{proof}

Theorem \ref{Comb} immediately follows from the combination of Lemma \ref{P2SP} and Lemma \ref{SPMP}.

\section{Realization of Weyl groups}
We retain the setting of \S 2. Let $R ( H )$ be the {\it complexified} representation ring of an algebraic group $H$. We have natural identification
$$Z \cong \{ ( g _1 B, g _2 B, X ) \in ( G / B ) ^2 \times \mathbb V : X \in g _1 \mathbb V ^+ \cap g _2 \mathbb V ^+ \}.$$
We have an inclusion $Z \subset F \times F$. Let $p _{ij} : F ^3 \rightarrow F ^2$ be the $(i, j)$-th projection ($1 \le i < j \le 3$). Since $p _{ij}$ is proper when restricted to $p _{12} ^{-1} ( Z ) \cap p _{23} ^{-1} ( Z )$, we have a well-defined convolution map
$$\star : K ( Z ) \otimes K ( Z ) \ni ( \mathcal F, \mathcal G ) \mapsto \sum _{i \ge 0} ( - 1 ) ^i [ \mathbb R ^i  ( p _{13} ) _* ( p _{12} ^* \mathcal F \otimes ^{\mathbb L} p _{23} ^* \mathcal G ) ] \in K ( Z ).$$

\begin{theorem}[Ginzburg]
The group $K ( Z )$ becomes an algebra via the convolution action. \hfill $\Box$
\end{theorem}

Let $\mathfrak t := \Lie T$. We define
$$\mathbb C [ \mathfrak t ] _0 := \mathbb C [ \mathfrak t ] / \left< \mathbb C [ \mathfrak t ] ^W _+ \right>,$$
where $\mathbb C [ \mathfrak t ] ^W _+ := \mathbb C [ \mathfrak t ] ^W \cap \mathfrak t ^* S ( \mathfrak t ^* )$. By the natural $W$-action on $\mathfrak t$, the space $\mathbb C [ \mathfrak t ] _0$ admits a $W$-action. Hence, we have their amalgamated product
$$\mathcal W := \mathbb C [ W ] \otimes \mathbb C [ \mathfrak t ] _0,$$
whose multiplication is given by $( w _1, f ) ( w _2, g ) := ( w _1 w _2, f w _1 ( g ))$.

For $\mathbb Z$-module $A$, we set $A _{\mathbb C} := \mathbb C \otimes _{\mathbb Z} A$.

\begin{theorem}\label{special}
We have an isomorphism $K ( Z ) _{\mathbb C} \cong \mathcal W$ as algebras.
\end{theorem}

\begin{proof}
Choose an element $(1, 1, -1, 1) \in G \times ( \mathbb C ^{\times} ) ^3$. It acts on $F _2 := G \times^{B} ( V _1 ^{\oplus 2} \oplus V _2 ) ^+$ with its fixed part isomorphic to $F \cong G \times ^B (  V _1 ^+ \oplus 0 \oplus V _2 ^+ ) \subset F _2$. By \cite{K} Corollary 2.13 and Remark 2.2 3), we conclude that $K ( Z ) _{\mathbb C}$ is isomorphic to the specialization of the three-parameter Hecke algebra of type $C _n ^{(1)}$ at $q _0 = - q _1 = q _2 = 1$ and $s = 1 \in G$. In particular, $K ( Z ) _{\mathbb C}$ is isomorphic to the quotient of the group ring of the affine Weyl group $W \ltimes X ^* ( T )$ by the ideal generated by the maximal ideal $\mathfrak m \subset R ( G )$ corresponding to $1 \in G$. (We regard $X ^* ( T ) \cong Q ^{\vee}$, where $Q ^{\vee}$ is the coroot lattice of $R$.) Here we have
$$R ( T ) / \mathfrak m R ( T ) = R ( T ) / \left< [ V ] - ( \dim V ) [ \mathbb C ] : V \in \mathrm{Rep} G  \right> \cong \mathbb C [ \mathfrak t ] ^W _0$$
as $W$-modules. Since $K ( Z ) _{\mathbb C} = \mathbb C [ W ] \otimes R ( T ) / \mathfrak m R ( T )$, the result follows.
\end{proof}

\begin{corollary}\label{inclmod}
We have a surjective map $\mathcal W \longrightarrow \!\!\!\!\! \rightarrow \mathbb C [ W ]$. In particular, we have an inclusion
$$\mathsf{Irr} W \subset \{ \text{simple } \mathcal W \text{-modules} \},$$
where the RHS is the set of isomorphism classes.
\end{corollary}

\begin{proof}
We retain the setting of the proof of Theorem \ref{special}. The maximal ideal $\mathfrak m _0 \subset \mathbb C [ T ]$ corresponding to $1 \in T$ is clearly $W$-invariant. We have $\mathbb C [ T ] \mathfrak m \subset \mathfrak m _0$. It follows that $\mathcal W / \mathcal W \mathfrak m _0 \cong \mathbb C [ W ]$ as desired.
\end{proof}

\begin{remark}
As is shown later (Proof of Theorem 6.5), the inclusion of Corollary \ref{inclmod} is in fact an equality. We can also deduce it directly from the structure of $\mathcal W$.
\end{remark}

\section{Rough classification of orbits}\label{sJNF}
We assume the setting of the previous sections. The goal of this section is to prove:

\begin{theorem}\label{surjMP}
The map
$$\mathbf J : \mathcal{MP} ( n ) \ni \vec{\lambda} \mapsto J ( \vec{\lambda} ) \in G \backslash \mathfrak N$$
is surjective.
\end{theorem}

\begin{proposition}[\cite{K} 1.9]\label{JNF}
Every element of $\mathfrak N$ is $G$-conjugate to a normal form. \hfill $\Box$
\end{proposition}

\begin{remark}
Theorem \ref{surjMP} does not follow from Proposition \ref{JNF} immediately since $\mathrm{Im} \mathbf J$ does not exhaust the set of normal forms modulo the $N _G ( T )$-action.
\end{remark}

Let $G _0 \subset G$ be the subgroup such that 1) $T \subset G _0$ and 2) the root system of $( G _0, T )$ is $\{ \pm ( \epsilon _i - \epsilon _j )\} _{i < j} \subset R$. We have $G _0 \cong \mathop{GL} ( n, \mathbb C )$. Let $X _{\epsilon _i - \epsilon _j} \in \mathrm{Lie} G _0 \subset \mathfrak g$ be a non-zero root vector of $T$-weight $\epsilon _i - \epsilon _j$ with standard normalization. We put $v _{i, j} := \exp ( X _{\epsilon _i - \epsilon _j} ) \in G _0$.

\begin{proof}[Proof of Theorem \ref{surjMP}]
Let $\underline{\lambda} = ( \lambda, a )$, where $\lambda$ is a partition of $n$ and $a = ( a _1, a _2, \ldots )$ is an integer sequence such that $0 \le a _p \le \lambda _p$ for each $p$. (We put $\vec{\lambda} = \underline{\lambda}$ if $\underline{\lambda}$ is a marked partition.) By Proposition \ref{JNF}, we assume
$$J = J ( \underline{\lambda} ) := \sum _i J ^{( a _i )} ( \lambda _i ) _{\lambda _i ^<} = J _1 + J ( \underline{\lambda} _0 ) \in \mathbb V ^+$$
with $J _1 \in V _1$ and $\underline{\lambda} _0 = ( \lambda, \mathbf 0)$. Let $i, j$ be two distinct integers such that $\lambda _i \ge \lambda _j$.
\begin{claim}
There exists elements $g ^+ _{i, j}, g ^- _{i, j} \in \mathrm{Stab} _{G _0} ( J ( \underline{\lambda} _0 ) )$ which induces unipotent transforms on $V _1$ so that
\begin{equation}
g ^+ _{i, j} \mathbf v [ \epsilon _{\lambda _j ^< + k} ] = \begin{cases} \mathbf v [ \epsilon _{\lambda _i ^< + k} ] + \mathbf v [ \epsilon _{\lambda _j ^< + k} ] & (0 < k \le \lambda _j ) \\ \mathbf v [ \epsilon _{\lambda _j ^< + k} ] & (otherwise) \end{cases} \label{gp}
\end{equation}
and
\begin{equation}
g ^- _{i, j} \mathbf v [ \epsilon _{\lambda _{i + 1} ^< - k} ] = \begin{cases} \mathbf v [ \epsilon _{\lambda _{i + 1} ^< - k} ] + \mathbf v [ \epsilon _{\lambda _{j + 1} ^< - k} ] & (0 \le k < \lambda _j ) \\ \mathbf v [ \epsilon _{\lambda _{i + 1} ^< - k} ] & (otherwise) \end{cases}.\label{gm}
\end{equation}
\end{claim}
\begin{proof}
It is straight-forward to see that $g ^{\pm} _{ij}$ satisfies (\ref{gp}) and (\ref{gm}). Define
$$g ^{+} _{ij} := \prod _{1 \le k \le \lambda _{j}} v _{\lambda ^< _i + k, \lambda ^< _j + k}, g ^{-} _{ij} := \prod _{0 \le k < \lambda _{j}} v _{\lambda ^< _{j + 1} - k, \lambda ^< _{i + 1} - k}\in G _0$$
fixes $J ( \underline{\lambda} _0 )$. By a weight comparison, these elements do not depend on the order of the product. The presentation of the transformations on $V _1$ follow immediately from this expression. We have
\begin{eqnarray*}
g ^+ _{i, j} J ( \underline{\lambda} _0 ) = J ( \underline{\lambda} _0 ) + \sum _{1 \le k < \lambda _j} X _{\epsilon _{\lambda _i ^< + k} - \epsilon _{\lambda _j ^< + k}} \mathbf v [ \epsilon _{\lambda _j ^< + k} - \epsilon _{\lambda _j ^< + k + 1} ]\\
+ \sum _{1 \le k < \lambda _j} X _{\epsilon _{\lambda _i ^< + k + 1} - \epsilon _{\lambda _j ^< + k + 1}} \mathbf v [ \epsilon _{\lambda _i ^< + k} - \epsilon _{\lambda _i ^< + k + 1} ] = J ( \underline{\lambda} _0 ).
\end{eqnarray*}
(Remember that the $G _0$-action on $\mathrm{Lie} G _0$ is the adjoint action.) This proves $g ^+ _{i, j} \in \mathrm{Stab} _{G _0} ( J ( \underline{\lambda} _0 ) )$. The case $g ^- _{i, j}$ follows from a similar equality.
\end{proof}
Assume that $J _1 = J _1 ^{\circ} + \mathbf v [ \epsilon _{\lambda _i ^{<} + a _i} ] + \mathbf v [ \epsilon _{\lambda _j ^{<} + a _j} ]$ with
$$\mathrm{supp} J _1 ^{\circ} \cap ( [ \lambda _i ^< + 1, \lambda _{i + 1} ^<] \cup [ \lambda _j ^< + 1, \lambda _{j + 1} ^<] ) = \emptyset.$$
We have
$$( g _{i, j} ^+ ) ^{-1} J = J - \mathbf v [ \epsilon _{\lambda _i ^{<} + a _j} ].$$
If $a _i < a _j$, then there exists an unipotent stabilizer $u \in \mathop{GL} ( \lambda _i ) \subset G _0$ of $J ( \underline{\lambda} _0 )$ such that
$$u ( g ^+ _{i, j} ) ^{-1} J = J ( \vec{\lambda} _0 ) + J _1 ^{\circ} - \mathbf v [ \epsilon _{\lambda _i ^{<} + a _j} ] - \mathbf v [ \epsilon _{\lambda _j ^{<} + a _j} ].$$
If $a _i = a _j$, then we have
$$( g ^+ _{i, j} ) ^{-1} J = J ( \vec{\lambda} _0 ) - \mathbf v [ \epsilon _{\lambda _j ^{<} + a _j} ].$$
In other words, we have $g ^+ _{i, j} u ( g ^+ _{i, j} ) ^{-1} J = J - \mathbf v [ \epsilon _{\lambda _{i} ^< + a _i} ]$ or $( g ^+ _{i, j} ) ^{-1} J = J - \mathbf v [ \epsilon _{\lambda _{i} ^< + a _i} ]$ when $a _i \le a _j$. Thus, if $a _i \le a _j$, then we have $G J ( \underline{\lambda} ) = G J ( \underline{\lambda} ^{\prime} )$ for $\underline{\lambda} ^{\prime} = ( \lambda, a ^{\prime} )$, where $a _p ^{\prime} = a _p$ ($p \neq i$) or $0$ ($p = i$).\\
By using $g ^- _{i, j}$ instead of $g ^+ _{i, j}$, we deduce:
If $\lambda _i - a _i \le \lambda _j - a _j$, then we have $G J ( \underline{\lambda} ) = G J ( \underline{\lambda} ^{\prime} )$ for $\underline{\lambda} ^{\prime} = ( \lambda, a ^{\prime} )$, where $a _p ^{\prime} = a _p$ ($p \neq j$) or $0$ ($p = j$).\\
We replace $( \lambda, a )$ by $( \lambda, a ^{\prime} )$ when one of the above two inequalities ($a _i \le a _j$ or $\lambda _i - a _i \le \lambda _j - a _j$) occur. Repeating these procedures for all possible pairs $( i, j )$ such that $\lambda _i \ge \lambda _j$, we obtain a marked partition $\vec{\mu} = ( \lambda, a ^{\prime \prime} )$ such that $G J ( \vec{\mu} ) = G J ( \underline{\mu} ) = G J ( \underline{\lambda} )$ as desired.
\end{proof}

\section{An exotic version of Springer correspondence}
We retain the setting of the previous section. The following result is not exactly the same as the original, but we can easily deduce it from the proof:

\begin{theorem}[Igusa \cite{Ig} Lemma 8]\label{Ig}
Let $\vec{\lambda} = ( \lambda, \mathbf 0) \in \mathcal P ( n ) \subset \mathcal{MP} ( n )$. Then, the reductive part of the stabilizer of $J ( \vec{\lambda} )$ is
$$L _{\lambda} := \mathop{Sp} ( 2n _1, \mathbb C ) \times \mathop{Sp} ( 2n _2, \mathbb C ) \times \cdots,$$
where the sequence $( n _1, n _2, \ldots )$ are the number of $\lambda _i$'s which share the same value. Moreover, we have
$$\mathsf{Res} ^G _{L _{\lambda}} V _1 = \bigoplus _{i \ge 1} V ( i ) ^{\oplus \lambda _i},$$
where $V ( i )$ is the external tensor product of a vector representation of $\mathop{Sp} ( 2 n _i )$ and trivial representations of $\mathop{Sp} ( 2 n _j )$ $(j \neq i)$. \hfill $\Box$
\end{theorem}

\begin{corollary}\label{stabconn}
For each $X \in \mathfrak N$, the $G$-stabilizer of $X$ is connected. 
\end{corollary}
\begin{proof}
We put $X = X _1 \oplus X _2 \in V _1 \oplus V _2$. Let $G _{1}$ and $G _{2}$ denote the $G$-stabilizers of $X _1$ and $X _2$. We show that $G _1 \cap G _2$ is connected. Let $G _2 = L _2 U _2$ be the Levi decomposition of $G _2$. The component group of $G _1 \cap G _2$ is the same as that of $G _1 \cap L _2$. By repeating the argument of \cite{K} Lemma 7.8, we conclude that $G _1 \cap L _2$ must be connected.
\end{proof}

\begin{theorem}[cf. \cite{CG} 8.9.3 and 8.4.8]\label{BBD}
There exists an algebraic stratification $\mathfrak O ^{\mu}$ of $\mathfrak N$ such that
$$\mu _* \mathbb C _{F} [ \dim F ] \cong \bigoplus _{\mathcal O \in \mathfrak O ^{\mu}} L _{\mathcal O} \boxtimes IC ( \mathcal O ),$$
where $L _{\mathcal O}$ is a vector space and $IC ( \mathcal O )$ is the minimal extension of $\mathbb C _{\mathcal O} [\dim \mathcal O]$.
\end{theorem}

\begin{proof}
By the BBD-G decomposition theorem (cf. Saito \cite{S} 5.4.8.2), we deduce that
$$\mu _* \mathbb C _{F} [ \dim F ] \cong \bigoplus _{d \in \mathbb Z, \mathcal O \in \mathfrak O ^{\mu}, \chi} L _{d, \mathcal O, \chi} \boxtimes IC ( \mathcal O, \chi ) [ d ] \in D ^b ( \mathfrak N ),$$
where $L _{d, \mathcal O, \chi}$ is a vector space, $\chi$ is a $G$-equivariant local system on $\mathcal O$, and $IC ( \mathcal O, \chi )$ is the minimal extension of $\chi$. 
The map $\mu$ is $G$-equivariant. It follows that $\mathfrak O _{\mathfrak N}$ is a refinement of $\mathfrak O ^{\mu}$. Since $\mu$ is semi-small with respect to $\mathfrak O _{\mathfrak N}$, it is also semi-small with respect to $\mathfrak O ^{\mu}$. Hence, we have $d \equiv 0$. (cf. \cite{CG} 8.9.2) Each strata $\mathcal O \in \mathfrak O ^{\nu}$ is smooth. In particular, the map $\pi _1 ( \mathbb O, *) \rightarrow \pi _1 ( \mathcal O, * )$ is surjective for the dense open $G$-orbit $\mathbb O \subset \mathcal O$. Therefore, we have $L _{d, \mathcal O, \chi} \neq 0$ only if $d = 0$ and $\chi = 1$ as desired.
\end{proof}

For $\mathbb O \in \mathfrak O _{\mathfrak N}$ and $X \in \mathbb O$, we define
$$M _{\mathbb O} := H _{\mathrm{codim} \mathcal O} ( \mu ^{- 1} ( X ), \mathbb C ) \text{ and } N _{\mathbb O} := \bigoplus _{m \ge 0} H _{m} ( \mu ^{- 1} ( X ), \mathbb C ),$$
where $\mathcal O \in \mathfrak O ^{\mu}$ is the strata such that $\mathbb O \subset \mathcal O$. By the Ginzburg theory, $N _{\mathbb O}$ is a $K ( Z )$-module. Since $G$ is connected, the $G$-conjugation of $X$ gives mutually isomorphic $K ( Z )$-modules. Thus, $N _{\mathbb O}$ is independent of the choice of $X$ as a $K ( Z )$-module with a grading.

\begin{theorem}[Chriss-Ginzburg]\label{cg}
Each $M _{\mathbb O}$ is a simple quotient of $N _{\mathbb O}$ as $K ( Z ) _{\mathbb C}$-modules if it is non-zero. Moreover, the set of isomorphism classes of non-zero modules in $\{ M _{\mathbb O} \} _{\mathbb O \in \mathfrak O _{\mathfrak N}}$ gives a complete collection of simple $K ( Z ) _{\mathbb C}$-modules.
\end{theorem}

\begin{proof}
Since $Z$ has a paving by affine spaces (see \cite{K} 1.5), it follows that $K ( Z )$ is spanned by algebraic cycles. By \cite{CG} 5.11.11, we have an isomorphism $K ( Z ) _{\mathbb C} \cong H _{\bullet} ( Z, \mathbb C )$ as convolution algebras. The first part follows from the combination of Theorem \ref{BBD}, \cite{CG} 8.9.8, and 8.9.14 (b). The second part follows by \cite{CG} 8.9.8.
\end{proof}

\begin{theorem}\label{main}
The assignment
$$\mathfrak O _{\mathfrak N} \ni \mathbb O \mapsto M _{\mathbb O} \in \mathsf{Irr} W$$
establish a one-to-one correspondence.
\end{theorem}

\begin{proof}
The subset $\mathfrak O _{\mathfrak N} ^{\circ} := \{ \mathbb O \in \mathfrak O _{\mathfrak N} : M _{\mathbb O} \neq 0 \}$ gives a map
$$\tau : \mathfrak O _{\mathfrak N} ^{\circ} \ni \mathbb O \mapsto M _{\mathbb O} \in \{\text{simple $K ( Z ) _{\mathbb C}$-modules} \}.$$
By Theorem \ref{cg}, this map must be surjective. Hence, we have
\begin{eqnarray*}
\# \{\text{simple $K ( Z ) _{\mathbb C}$-modules} \} \le \# \mathfrak O _{\mathfrak N} ^{\circ} \le \# \mathfrak O _{\mathfrak N} \le \# \mathcal{MP} ( n )\\ = \# \mathcal P _2 ( n ) = \# \mathsf{Irr} W \le \# \{\text{simple $K ( Z ) _{\mathbb C}$-modules} \}
\end{eqnarray*}
by Theorem \ref{surjMP}, Theorem \ref{Comb}, and Corollary \ref{inclmod}. This implies that all the inequalities are in fact an equality. Therefore, the map $\tau$ is defined at the whole of $\mathfrak O _{\mathfrak N}$ and the map is injective. Since every simple $W$-module give rise to a simple $K ( Z ) _{\mathbb C}$-module, the result follows.
\end{proof}

\begin{corollary}\label{cormain}
The set $\mathfrak O _{\mathfrak N}$ is in one-to-one correspondence with $\mathcal P _2 ( n )$. \hfill $\Box$ 
\end{corollary}

\section{Regularity conditions of parameters}
This section might also be viewed as an continuation of \cite{K} with the knowledge of this paper. In the below, we use notation of \cite{K} freely only by indicating pointers to them. We assume the setting of \S 4. Let $a = ( s, q _1, q _2 ) \in T \times ( \mathbb C ^{\times} ) ^2$. (This is equivalent to assume $a \in G \times ( \mathbb C ^{\times} ) ^2$ is semi-simple by taking an appropriate conjugate.)

Consider the following condition:
\begin{itemize}
\item[$(\sharp) _1$] $q _2$ is not a root of unity of order $\le 2n$;
\item[$(\sharp) _2$] Let $\Psi ( s )$ be the set of $s$-eigenvalues in $V _1$. For each $c _1, c _2 \in \Psi ( s )$, we have $c _1 = q _2 ^m c _2 ^{\pm 1}$ for some integer $m$ and
$$c _1, q _2 ^{-1} c _1, \ldots, q _2 ^{-m} c _1 \in \Psi ( s ).$$
\end{itemize}

\begin{theorem}\label{oneps}
Assume $(\sharp)$. There exists a one-parameter subgroup $\psi : \mathbb C ^{\times} \to T \times ( \mathbb C ^{\times} ) ^2$ such that:
$$Z _{G \times ( \mathbb C ^{\times} ) ^2} ( \psi ( r ) ) = Z _{G \times ( \mathbb C ^{\times} ) ^2} ( a ), \text{ and } \mathfrak N ^{\psi ( r )} = \mathfrak N ^a$$
for a generic choice of $r \in \mathbb R _{>0}$.
\end{theorem}
\begin{proof}
By \cite{K} \S 4, the condition $(\sharp)$ implies that the setting is governed by the relations and values of $q _1 = e ^{r _1}, q _2 = e ^{r _2},$ and $s = \exp ( \sum _{i = 1} ^n \lambda _i \epsilon _i )$. In particular, we can rearrange their values to be $r _i, \lambda _i \in \mathbb R$ without changing the $\mathfrak N ^a$ and $F ^a$ from the original ones. 
\end{proof}

In the setting of Theorem \ref{oneps}, we define $A = A ( r )$ to be the Zariski closure of $\{\psi ( m r ) : m \in \mathbb Z\} \subset T \times ( \mathbb C ^{\times} ) ^2$.

\begin{corollary}\label{atorus}
Keep the setting of Theorem \ref{oneps}. For a generic choice of $r \in \mathbb R _{>0}$,  the torus $A ( r )$ is connected.
\end{corollary}

\begin{proof}
If $\psi ( r )$ is not sitting in the identity component of $A$, then so does each of $\psi ( r / m )$ ($m \in \mathbb Z _{>0}$). This is impossible since $A$ has only finitely many connected components by definition. This contradiction implies $\psi ( r ) \in A$, which in turn yields that $A$ is connected.
\end{proof}

We assume $(\sharp)$ and the setting of Theorem \ref{oneps} until Theorem \ref{reg}.

For each $m \ge 0$, let $EA _m := ( \mathbb C ^m \backslash \{ 0 \} ) ^{\dim A}$ be a variety such that $i$-th $\mathbb C ^{\times}$-factor of $A = ( \mathbb C ^{\times} ) ^{\dim A}$ acts as dilation of the $i$-th factor for each $1 \le i \le m$. By the standard embedding $\mathbb C ^m \hookrightarrow \mathbb C ^{m + 1}$ sending $( x )$ to $( x , 0 )$, we form a sequence of $A$-varieties
$$\emptyset = EA _0 \hookrightarrow EA _1 \hookrightarrow EA _2 \hookrightarrow \cdots.$$
We define $EA := \varinjlim _m EA _m$, which is an ind-quasiaffine scheme with free $A$-action. Since $EA$ is contractible manifold with respect to the classical topology, we regard $EA$ as the classifying space of $A$. For a $A$-variety $\mathcal X$, we set
$$\mathcal X _A := \triangle A \backslash \left( EA \times \mathcal X \right).$$
We have a forgetful map
$$f ^A _{\mathcal X} : \mathcal X _A \rightarrow BA = A \backslash EA.$$
Let $\mathbb D _X ^A$ be the relative dualizing sheaf with respect to $f _{\mathcal X}$. We define
$$H _{i} ^A ( \mathcal X ) \cong H ^{-i} ( \mathcal X _A, \mathbb D _X ^A ).$$
We have the Leray spectral sequence
$$H ^i ( BA ) \otimes H _{j} ( \mathcal X ) \Rightarrow H _{- i + j} ^A ( \mathcal X ).$$
In the below, we understand that $H _{\bullet} ^A ( \mathcal X ) := \bigoplus _{m} H _{m} ^A ( \mathcal X )$. The projection maps $p _i : Z _A \rightarrow F _A$ ($i = 1, 2$) equip $H _{\bullet} ^A ( Z )$ a structure of convolution algebra. It is straight-forward to see that the diagonal subsets $\triangle F \subset Z$ and $( \triangle F ) _A \subset Z _A$ represents $1 \in H _{\bullet} ( Z )$ and $1 \in H _{\bullet} ^A ( Z )$, respectively.

\begin{lemma}
The algebra $H _{\bullet} ^A ( Z )$ contains $H _{\bullet} ( Z )$ as its subalgebra. In particular, we have $\mathbb C [ W ] \subset H _{\bullet} ^A ( Z )$ as subalgebras. Moreover, the center of $H _{\bullet} ^A ( Z )$ contains $H ^{\bullet} ( BA ) [ ( \triangle F ) _{A} ] \subset H _{\bullet} ( Z )$.
\end{lemma}

\begin{proof}
In the Leray spectral sequence
$$H ^i ( BA ) \otimes H _{j} ( Z ) \Rightarrow H _{- i + j} ^A ( Z ),$$
we have $H ^{odd} ( BA ) = 0$ (by Corollary \ref{atorus}) and $H _{odd} ( Z ) = 0$ (since $Z$ is paved by affine spaces). It follows that this spectral sequence degenerates at the level of $E _2$-terms. Moreover, the image of the natural map $\imath : H _j ( Z ) \hookrightarrow H _{j} ^A ( Z )$ represents cycles which are locally constant fibration over the base $BA$. It follows that the map $\imath$ is an embedding of convolution algebras.\\
Multiplying $H ^{\bullet} ( BA )$ is an operation along the base $BA$, which commutes with the convolution operation (along the fibers of $f ^A _{Z}$). It follows that $H ^{\bullet} ( BA ) \rightarrow H ^{\bullet} ( BA ) [ ( \triangle F ) _A ] \subset H _{\bullet} ^A ( Z )$ is central subalgebra as desired.
\end{proof}

The following result is a consequence of Borho-MacPherson's argument applied to the sheaf $\mu _* \mathbb C$ as in the previous section.

\begin{theorem}[Borho-MacPherson cf. \cite{CG} \S 8.8]\label{BorMac}
Let $y, y ^{\prime} \in \mathfrak N$. Then, we have
$$\hskip 20mm [ H _{\bullet} ( \mu ^{-1} ( y ) ) : M _{G y ^{\prime}} ] = \begin{cases} 0 & (y \not\in \overline{G y ^{\prime}}) \\ 1 & (y \in G y ^{\prime}) \end{cases}$$
as $W$-modules. \hfill $\Box$
\end{theorem}

\begin{lemma}\label{LerayM}
Let $y \in \mathfrak N$. The Leray spectral sequence
$$H ^{\bullet} ( BA ) \otimes H _{\bullet} ( \mu ^{-1} ( y ) ) \longrightarrow H _{\bullet} ^A ( \mu ^{-1} ( y ) )$$
induces a map of $H _{\bullet} ^A ( Z )$-module by letting $H ^{\bullet} ( BA )$ act only on the first term of the LHS and letting $H _{\bullet} ( Z )$ act only on the second term of the LHS.
\end{lemma}

\begin{proof}
The LHS is the cohomology space of a sheaf of $H _{\bullet} ( Z )$-modules on $BA$. For any contractible set $\mathcal U \subset BA$, we have
$$H _{\bullet} ^A ( \mu ^{-1} ( y ) ) \longrightarrow H _{\bullet} ( ( f ^A _{\mu ^{-1} ( y )} ) ^{-1} ( \mathcal U ) ) \cong H _{\bullet} ( \mu ^{-1} ( y ) )$$
as $H _{\bullet} ( Z )$-modules. It follows that the composition map
$$H ^0 ( BA ) \otimes H _{\bullet} ( \mu ^{-1} ( y ) ) \rightarrow H _{\bullet} ^A ( \mu ^{-1} ( y ) ) \rightarrow H _{\bullet} ( \mu ^{-1} ( y ) )$$
is a $H _{\bullet} ( Z )$-module map. Since the map
$$H ^{\bullet} ( BA ) \otimes H _{\bullet} ( \mu ^{-1} ( y ) ) \longrightarrow H _{\bullet} ^A ( \mu ^{-1} ( y ) )$$
is a $H ^{\bullet} ( BA )$-module map, we conclude that it is a $H ^A _{\bullet} ( Z )$-module map as desired. 
\end{proof}

Let $y \in \mathfrak N ^a$. We put
$$\nabla ^A _{(a, y)} := H ^A _{\bullet} ( \mu ^{-1} ( y ) ) / H ^{\bullet} ( BA ) \text{-torsion}.$$

\begin{lemma}
The $H ^A _{\bullet} ( Z )$-module $\nabla ^A _{(a, y)}$ contains $W$-module $H ^{\bullet} ( BA ) \otimes M _{Gy}$.
\end{lemma}

\begin{proof}
By Borho-MacPherson's theorem, the $W$-module $M _{Gy}$ appears in $H _{\bullet} ( \mu ^{-1} ( y ) )$ with multiplicity one. It follows that the $M _{Gy}$-isotypical component of the Leray spectral sequence
$$H ^{\bullet} ( BA ) \otimes H _{\bullet} ( \mu ^{-1} ( y ) ) \longrightarrow H _{\bullet} ^A ( \mu ^{-1} ( y ) )$$
is $E _2$-degenerate. This implies that the $M _{Gy}$-isotypical component of $H _{\bullet} ^A ( \mu ^{-1} ( y ) )$ is a free $H ^{\bullet} ( BA )$-module, which cannot be torsion.
\end{proof}

We put $\mathfrak a := \mathrm{Lie} A$. We have $H ^{\bullet} ( BA ) \cong \mathbb C [ \mathfrak a ]$. By inverting all monomials which are perpendicular to $a$, we obtain a localized algebra $H ^{\bullet} ( BA ) _a := \mathbb C [ \mathfrak a ] _a$. For a $A$-variety $\mathcal X$, we put
$$H ^A _{\bullet} ( \mathcal X ) _a := H ^{\bullet} ( BA ) _a \otimes _{H ^{\bullet} ( BA )} H ^A _{\bullet} ( \mathcal X ).$$
Let $H ^+ ( BA ) := \bigoplus _{m > 0} H ^m ( BA )$.

\begin{proposition}\label{degHa}
We have an isomorphism of convolution algebras:
\begin{eqnarray}
H ^A _{\bullet} ( Z ^a ) _a \cong H ^A _{\bullet} ( Z ) _{a}.\label{localH}
\end{eqnarray}
Moreover, the quotient space
$$\nabla ^{\prime} _{(a, y)} := \nabla ^{A} _{(a, y)} / H ^+ ( BA ) \nabla ^{A} _{(a, y)}$$
admits a $H _{\bullet} ( Z ^a )$-module structure such that
\begin{enumerate}
\item $\nabla ^{\prime} _{(a, y)}$ is a subquotient of $H _{\bullet} ( \mu ^{-1} ( y ) ^a )$ as $H _{\bullet} ( Z ^a )$-modules;
\item The map $\mathbb C [ W ] \subset H _{\bullet} ( Z ) \hookrightarrow H ^A _{\bullet} ( Z ^a ) _a$ defines a $W$-module structure on $\nabla ^{\prime} _{(a, y)}$;
\item $\nabla ^{\prime} _{(a, y)}$ is a quotient module of $H _{\bullet} ( \mu ^{-1} ( y ) )$ as $W$-modules.
\end{enumerate}
\end{proposition}

\begin{proof}
Let $R ( A ) _a$ be the localization of $R ( A )$ at the point $a$. By the Thomason localization theorem (see e.g. \cite{CG} \S 8.2), we have an isomorphism
$$R ( A ) _a \otimes _{R ( A )} K ^A ( Z ^a ) \cong R ( A ) _a \otimes _{R ( A )} K ^A ( Z )$$
as algebras. For each of $\mathcal X = Z$, or $Z ^a$, we have a dense open embedding
$$K ^A ( \mathcal X ) \hookrightarrow \varprojlim _m K ^A ( EA _m \times \mathcal X ) \cong \varprojlim _m K ( A \backslash ( EA _m \times \mathcal X ) ).$$
We regard the RHS as a substitute of $K ( \mathcal X _A )$. It follows that the Chern character map relative to $BA$ gives an isomorphism
$$\mathbb C [[ \mathfrak a ]] _{a} \otimes _{\mathbb C [ \mathfrak a ] _a} H ^A _{\bullet} ( Z ^a ) _a \cong \mathbb C [[ \mathfrak a ]] _{a} \otimes _{\mathbb C [ \mathfrak a ] _a} H ^A _{\bullet} ( Z ) _{a},$$
where $\mathbb C [[ \mathfrak a ]] _{a}$ is the formal power series ring of $\mathbb C [ \mathfrak a ]$ along $a$. By restricting this to the sum of vectors of finitely many degrees, we obtain (\ref{localH}). The second assertion is automatic by letting $H _{\bullet} ( Z ^a )$ act by
$$H _{\bullet} ( Z ^a ) \longrightarrow H _{\bullet} ^A ( Z ^a ) \longrightarrow H _{\bullet} ^A ( Z ^a ) / H ^+ ( BA ) H _{\bullet} ^A ( Z ^a ).$$
By a similar argument using Lemma \ref{LerayM}, we deduce that $\nabla ^{\prime} _{(a, y)}$ is a subquotient of $H _{\bullet} ( \mu ^{-1} ( y ) ^a )$, which implies 1). Now we verify 2). Since (\ref{localH}) is an algebra isomorphism, it follows that $1 \in \mathbb C [ W ]$ goes to $1 \in H _{\bullet} ( Z ^a )$. It follows that each of $s _i$ goes to a non-zero element of $H _{\bullet} ( Z ^a )$ with its square equal to $1$. By construction, there exists $f _i \in \mathbb C ( \mathfrak a )$ ($i = 1, \ldots n$) such that $1, f _1 s _1, \ldots f _n s _n \in H _{\bullet} ^A ( Z ^a ) _a$ define linearly independent vectors in $H _{\bullet} ( Z ^a )$. It follows that $f _i ^2 \in \mathbb C [ \mathfrak a ]$. This forces $f _i \in \mathbb C [ \mathfrak a ]$, which implies that the images of $1, s _1, \ldots, s _n \in H _{\bullet} ( Z ^a )$ are linearly independent. This verifies 2). Since $\mathbb C [W]$ is a semi-simple algebra, we have $W$-module morphisms whose composition is surjective
$$H _{\bullet} ( \mu ^{-1} ( y ) ) \hookrightarrow H ^{\bullet} ( BA ) \otimes H _{\bullet} ( \mu ^{-1} ( y ) ) \longrightarrow \nabla ^{\prime} _{(a, y)}.$$
This verifies 3) as desired.
\end{proof}

\begin{proposition}\label{exclorbits}
Let $\mathcal O \subset \mathfrak N$ be a $G$-orbit. For any two distinct $G ( s )$-orbits $\mathcal O _1, \mathcal O _2 \subset \mathcal O ^a$, we have
$$\overline{\mathcal O _1} \cap \mathcal O _2 = \emptyset.$$
\end{proposition}

\begin{proof}
By the description of $G$-orbits of $\mathfrak N$, we deduce that the scalar multiplication of a normal form of $\mathfrak N$ is achieved by the action of $T$. It follows that each $G ( s )$-orbit of $\mathfrak N ^a$ is a $Z _{G \times ( \mathbb C ^{\times} ) ^2} ( a )$-orbit. Let $y \in \mathcal O _1$. Let $G _y$ be the stabilizer of $y$ in $G \times ( \mathbb C ^{\times} ) ^2$. Assume that $\mathcal O _2 \cap \overline{\mathcal O _1} \neq \emptyset$ to deduce contradiction. Since $\mathcal O _2$ is a $Z _{G \times ( \mathbb C ^{\times} ) ^2} ( a )$-orbit, we have $\mathcal O _2 \subset \overline{\mathcal O _1}$. Fix $y _2 \in \mathcal O _2$. Consider an open neighborhood $\mathcal U$ of $1$ in $G$ (as complex analytic manifolds). Then, $\mathcal U y _2 \in \mathcal O$ is an open neighborhood of $y _2$. It follows that $\mathcal U y _2 \cap \mathcal O _1 \neq \emptyset$. We put $\mathfrak g _{a, y _2} := \mathrm{Lie} G _{y _2} + \mathrm{Lie} Z _{G \times ( \mathbb C ^{\times} ) ^2} ( a )$. We have
$$N _{\mathcal O _2 / \mathcal O, y _2} = \mathfrak g / \mathfrak g _{a, y _2}.$$
Every non-zero vectors of $N _{\mathcal O _2 / \mathcal O, y _2}$ is expressed as a linear combination of eigenvectors with respect to the $a$-action. These $a$-eigenvectors can be taken to have non-zero weights and does not contained in $G _{y _2}$. It follows that
$$\mathcal U y _2 \cap \mathcal O _1 \not\subset \mathbb V ^a,$$
which is contradiction (for an arbitrary sufficiently small $\mathcal U$). Hence, we have necessarily $\mathcal O _2 \cap \overline{\mathcal O _1} = \emptyset$ as desired.
\end{proof}

\begin{theorem}\label{reg}
We have the following one-to-one correspondence
$$G ( s ) \backslash \mathfrak N ^a \ni G ( s ) y \mapsto L _{(a, y)} \in \mathsf{Irr} H _{\bullet} ( Z ^a ),$$
where $L _{(a, y)}$ is the unique $H _{\bullet} ( Z ^a )$-irreducible constituent of $\nabla _{(a, y)}$ which contains $W$-module $M _{Gy}$.
\end{theorem}

\begin{proof}
By \cite{K} Theorem 1.8, the number of $G ( s )$-orbits in $\mathfrak N ^a$ is finite. By \cite{K} Theorem 8.1, each $G ( s )$-orbit of $\mathfrak N ^a$ corresponds to at most one irreducible $H _{\bullet} ( Z ^a )$-module. We prove the assertion by the induction on the closure relation of the orbits. By Ginzburg's theory \cite{CG} \S 8.4 and Proposition \ref{degHa} 1), we deduce that $L _{(a, y)}$ does not contain a $W$-module which do not appear in $\nabla ^{\prime} _{(a, y ^{\prime})}$ for some $y \in \overline{G ( s ) y ^{\prime}}$. By Lemma \ref{exclorbits} and Theorem \ref{BorMac}, we deduce that $\nabla ^{\prime} _{(a, y)}$ carries a $W$-module $M _{G y}$, which does not contained in $\nabla ^{\prime} _{(a, y ^{\prime})}$ for every $y ^{\prime} \in \mathfrak N ^a$ such that $y \in \overline{G y ^{\prime}} \backslash G y ^{\prime}$. In particular, we conclude $L _{(a, y)} \neq 0$ as desired.
\end{proof}

We forget the assumption $(\sharp)$.

Recall that an extended affine Hecke algebras of type $B _n$ with two-parameters $( - q _1, q _1, q _2 )$ is the quotient of an affine Hecke algebra $\mathbb H$ of type $C _n$ with three parameters $( q _0, q _1, q _2 )$ by the two-sided ideal $( q _0 + q _1 )$. (See \cite{K} 2.1--2.2 for more detail.)

\begin{corollary}\label{BDL}
For an extended affine Hecke algebras of type $B _n$ with two-parameters $( - q _1, q _1, q _2 )$, the regularity condition of parameters holds automatically unless $- q _1 ^2 \neq q _2 ^{\pm m}$ for $0 \le m < n$ or $q _2 ^l \neq 1$ for $1 \le l < 2n$.
\end{corollary}

\begin{proof}
Applying \cite{K} Corollay 3.10, we can assume $(\sharp)$ freely. Hence, if we have either $V _1 ^{( s, q _1 )} = \{ 0 \}$ or $V _1 ^{( s, - q _1 )} = \{ 0 \}$, then Theorem \ref{reg} implies the result. Here $s$ has at most $2n$-eigenvalues and every two eigenvalues $\xi _1, \xi _2$ are connected by $\xi _1 = q _2 ^{m} \xi _2 ^{\pm 1}$ by some $- n < m < n$ by $(\sharp) _2$. Therefore, we cannot have both $V _1 ^{( s, q _1 )} \neq \{ 0 \}$ and $V _1 ^{( s, - q _1 )} \neq \{ 0 \}$ simultaneously by $- q _1 ^2 \neq q _2 ^{\pm m}$.
\end{proof}

{\bf Acknowledgement:} The author would like to thank Professor Hisayosi Matumoto for discussion on this topic. The author also likes to express his thanks to Professors Michel Brion, Toshiyuki Tanisaki and Takuya Ohta for their advice.

\end{document}